\input amstex
\documentstyle{amsppt}
\magnification=\magstep1

\def\thetitle{Remarks on certain composita of fields}
\def\thedate{13 December 2011}

\let\endignore=\relax
\long\def\ignore#1\endignore{\relax}
%
%

\def\Capelli    {1}
\def\Douady     {2}   
\def\Harbater   {3}   
\def\HaranJarden{4}   
\def\Lorenz     {5}   
\def\Pop        {6}   
\def\Serre      {7}   
\def\Witt       {8}   
\def\ZS         {9}

\define\ZZ{{\Bbb Z}}

\define\FF{{\Bbb F}}

\TagsOnRight

\def\Soc{\operatorname{Soc}}
\def\Gal{\operatorname{Gal}}
\def\Cycl{\operatorname{C}}
 \def\wpinv#1{\wp^{-1}#1}
\define\proot#1{\root p\of{#1}}
\def\Foverline(#1){\hskip2,5pt \widetilde{\hskip-2.5pt 
       F(#1)\hskip-1pt }\hskip2pt\relax}
\def\Coverline(#1){\hskip2,5pt \widetilde{\hskip-2.5pt 
       C(#1)\hskip-1pt }\hskip2pt\relax}
\def\cupeq{\cup\hbox{\kern 1pt \vrule height6.3pt \relax}}

\topmatter
   \date \thedate \enddate
 \title  \thetitle  \endtitle
 \author Christian U\. Jensen and Anders Thorup
\endauthor
\affil 
Department of Mathematical Sciences, University of Copenhagen
\endaffil 
\address  Universitetsparken 5, DK--2100 Copenhagen, Denmark 
\endaddress 
\email cujensen\@math.ku.dk, \ thorup\@math.ku.dk \endemail
\catcode`\@=11
   \def\subjyear@{2010}%
\catcode`\@=\active
 \abstract\nofrills{\bf Abstract.}\ \ 
 Let $L$ and $M$ be two algebraically closed fields contained in some
common larger field.  It is obvious that the intersection $C=L\cap M$
is also algebraically closed.  Although the compositum $ L M$ is
obviously perfect, there is no reason why it should be algebraically
closed except when one of the two fields is contained in the other. We
prove that if the two fields are strictly larger that $C$, and
linearly disjoint over $C$, then the compositum $LM$ is not
algebraically closed; in fact we shall prove that the Galois group of the
maximal abelian extension of $LM$ is the free pro-abelian group of
rank $|LM|$, and that the free pro-nilpotent group of rank $|C |$ 
can be realized as a Galois group over $L M$.

The above results may be considered as the main contribution of this
article but we obtain some additional results on field composita that
might be of independent interest. 
 \endabstract
\endtopmatter

\document

\subhead  Introduction
\endsubhead

 Let $L$ and $M$ be two algebraically closed fields contained in some
common larger field.  Clearly, the intersection $C=L\cap M$
is also algebraically closed.  Offhand there is no reason why the
compositum $LM$ should be algebraically closed unless one the fields
contains the fields is contained in the other.
In this paper we show that if $L$ and $M$ are proper extensions of
$C$, and linearly disjoint over $C$, then $LM$ is not algebraically
closed (not even separably closed).  In particular, if one of the
fields $L$ or $M$ has transcendence degree $1$ over $C$, then $LM$ is
not algebraically closed.

The three first sections are basically elementary.  Section 1 states
as the main observation (Theorem 1.3) that the class of algebraically
relatively closed extensions is stable under formation of composita
with separably generated extensions.  This result plays an essential
role in the paper. 

Section 2 presents a detailed analysis of $p$-socles of field
extensions ($p$ is a given prime) and their relationship to
$p$-Frattini groups. In addition the section contains a short survey of
Kummer theory and Artin--Schreier theory.

Section 3 combines sections 1 and 2 to show that elementary abelian
$p$-groups up to a certain rank can be realized by an
explicit construction  as Galois groups over $LM$ when $L/C$
and $M/C$ are proper extensions and linearly disjoint over $C$.

 Section 4 contains the two main results. For the first one, we use the
results on elementary abelian $p$-extensions to describe the maximal
abelian extension of $LM$: Its Galois group is the free pro-abelian
group of rank $|LM|$. 

Finally, for the second main result we use the rather deep
Douady--Harbater--Pop theorem on the absolute Galois group of a
function field in one variable, and a classical theorem of Witt on
$p$-extensions of a field of characteristic $p$.  Combined with the
results of Section 3 we show for any prime $p$ 
that the free pro-$p$ group of rank $|C|$ is realizable as a
Galois group over $LM$.
 As a consequence, the free pro-nilpotent group of rank $|C|$ is
realizable as a Galois group over $LM$.

 The precise results, Theorems 4.1 and 4.2
may be
considered as the main contribution of this article, but some of the
results on field composita in Sections 1 and 3 may be of independent
interest.

\subhead
Section 1. Auxiliary results on relative algebraic closures of fields
\endsubhead

 For the convenience of the reader we recall in this section some
concepts and state some elementary results from classical
field theory.

\definition{Definitions 1.1} A subfield $F$ of the
field $L$ is called {\it algebraically closed relative to} $L$ if any
element of $L$ that is algebraic over $F$ lies in $F$, i.e., if $F$
equals the algebraic closure of $F$ in $L$.

 A field extension $M$ of $F$ is called
{\it separably generated\/} if there exists a transcendency basis $T$
of $M$ with respect to $F$ such that $M$ is an algebraic separable
extension of $F(T)$.  

 Recall that two extensions $L_1/F$ and $L_2/F$ are called {\it
linearly disjoint\/} (over $F$) if elements of $L_2$ that are linearly
independent over $F$ are also linearly independent over $L_1$.  More
symmetrically, the two extensions are linearly disjoint when for any
two families $(x_i)$ in $L_1$ and $(y_j)$ in $L_2$, both linearly
independent over $F$, the family of all products $(x_iy_j)$ is
linearly independent over $F$.  
 In particular, if $T$ and $U$ are two sets of independent
variables over $F$ then $F(T)$ and $F(U)$ are linearly disjoint over
$F$ if and only if the union $S\cup T$ is a set of independent
variables over $F$. 
 \enddefinition

\remark{Remark 1.2} Two finite extensions $L_1/F$ and $L_2/F$ are
linearly disjoint if and only if
$\dim_F L_1L_2=(\dim _FL_1)(\dim_FL_2)$. As a consequence, if the
extensions are subfields of a finite Galois extension of $F$ with
Galois group $G$, and $H_1$ and $H_2$ are the subgroups of
automorphisms in $G$ fixing, respectively, $L_1$ and $L_2$, then the
extensions are linearly disjoint if and only if $G=H_1H_2$, where
$H_1H_2$ is the subset of $G$ consisting of products $h_1h_2$ with
$h_1\in H_1$ and $h_2\in H_2$. 

 Finally, assume that one of the two extensions, say $L_1/F$, is a
Galois extension with group $G:=\Gal(L_1/F)$. We will make extensive
use of the Translation Theorem of Galois theory, see
\cite{\Lorenz, Chapter 12, Theorem 1, p.  115}: The extension
$L_1L_2/L_2$ is Galois, and restriction of automorphisms defines an
injection $\Gal(L_1L_2/L_2)\hookrightarrow \Gal (L_1/F)$.  Moreover,
the two extensions $L_1/F$ and $L_2/F$ are linearly disjoint iff
$L_1\cap L_2=F$, iff the injection of Galois groups is an isomorphism.
 \endremark

The main result in this section is the following.

\proclaim {Theorem 1.3}
 Let $M$ and $L$ be fields contained in some
larger field and assume that $L$ and $M$ are linearly disjoint over
their intersection $F = M\cap L$.
 If M is a separably generated over $F$ and $F$ is algebraically
closed relative to $L$ then $M$ is algebraically closed relative to 
the compositum $M L$. 
\endproclaim

 Clearly, to prove the Theorem it suffices to prove it when $M/F$ is a
purely transcendental extension, and when $M/F$ is an algebraic
separable extension.  Moreover, it suffices to prove it when $M/L$ is
finitely generated, and we may in fact assume that $M/F$ is generated
by a single element. We separate the proof of the Theorem into the
transcendental and the algebraic part.

\proclaim{Proposition 1.4} If $F$ is a subfield of $L$ and $F$ is
algebraically closed relative to $L$, then the field of rational functions
$F(T)$, where $T$ is set of algebraic independent elements over $L$,
is algebraically closed relative to 
the field of rational functions $L(T)$. 
\endproclaim 
\demo{Proof}
As noted above we may assume the $T$ consist of a single
element $t$, transcendental over $L$.
Let $\varphi\in L(t)$, $\varphi\neq 0$, be algebraic over $F(t)$, so
there is an equation 
 $$ f_n\varphi^n + \cdots + f_0 = 0  \tag**$$
 where $f_i\in F(t)$ and $f_n \neq 0$. We must prove that $\varphi\in
F(t)$.
 
Multiplying the functions $f_i$ by a common denominator in $F[t]$, we
may assume that each $f_i$ in (*) is a polynomial in $F[t]$. We
multiply (*) by $f_n^{n-1}$ and obtain an equation of degree $n$
showing that $f_n\varphi $ is integral over $F[t]$. Since $\varphi $
lies in $F(t)$ if $f_n\varphi $ does, we may assume that $\varphi $ is
integral over $F[t]$. In particular, $\varphi $ is integral over
$L[t]$, hence $\varphi $ must lie in $L[t]$, because $L[t]$ being UFD
is integrally closed in its quotient field.  \smallskip Next we show
that the coefficients of $\varphi $ are algebraic over $F$. From the
assumption that $F$ is algebraically closed relative to $L$ it then
follows that $\varphi $ lies in $F[t]$.  \smallskip Let $\alpha$ be
the leading coefficient in $\varphi $. By induction on the degree of
$\varphi $ it suffices to show that $\alpha$ is algebraic over $F$.
Let $N$ be highest degree of the polynomials $f_i\varphi ^i$ for $ i =
0,1\dots,n$. Using that the coefficient of $t^N$ on the left hand side
of (*) is 0 we obtain
 $$ a_n\alpha^n +\cdots + a_0 = 0 
 $$
 where $a_i$ is the leading coefficient of $f_i$ if $f_i\varphi ^i$
has degree $N$ and $a_i = 0$ otherwise. By the choice of $N$ at least
one $a_i$ is $\neq 0$. Thus $\alpha$ is algebraic over $F$. \hfill
$\square$ 
 \enddemo

\proclaim{Proposition 1.5} If $F$ is a subfield of $L$ and $F$ is
algebraically closed relative to $L$ and $M$ is a finite separable
extension of $F$, then $M$ is algebraically closed relative to $M L$.
 \endproclaim
 \demo{Proof} Since $M$ is a finite separable extension of $F$ there
exists $\theta\in M$ such that $M = F(\theta)$. We have to show that
$F(\theta)$ is algebraically closed relative to $L(\theta)$.

Let $f\in F[X]$ be the minimal polynomial of $\theta$ over $F$, of
degree $n$ say.  Since $F$ is algebraically closed relative to $L$, it
follows that $f$ is irreducible over $L$. Thus $[L(\theta):L] =
[F(\theta):F]=n$, and any $\lambda\in L(\theta)$ be can be written 
 $$ \lambda = \ell_0 + \ell_1\theta+\cdots +
\ell_{n-1}\theta^{n-1}.\tag**
 $$
 where $\ell_0,\ell_1,\dots,\ell_{n-1}$ lie in $L$.

 Assume that $\lambda$ is algebraic over $F(\theta)$.  We must show
that $\ell_0,\ell_1,\dots,\ell_{n-1}$ actually lie in $F$.

We can write $f(X) = \prod_{i =1}^n(X - \theta_i)$, where $\theta_i$,
$1\leqq i\leqq n$ (by the separability assumption) are distinct
elements in the Galois closure $M'$ of $M$.  We may assume that
$\theta = \theta_1$.  For each $i$, $1\leqq i\leqq n$, there exists an
automorphism $\sigma_i\in \Gal(M'/F)$ such that $\sigma_i(\theta) =
\theta_i$.  Since $M'\cap L=F$, the automorphism $\sigma
_i\in \Gal(M'/F)$ extends uniquely to an automorphism (also denoted
$\sigma _i$) in $\Gal(M'L/L)$, see Remark 1.2,

Clearly for $1\leqq i\leqq n$ we have 
 by (**)
 $$ \sigma_i(\lambda) = \ell_0 + \ell_1\theta_i+\cdots +
\ell_{n-1}\theta_i{}^{\!n-1}.\tag***
 $$
 Since $\lambda$ is algebraic over $F(\theta)$ and $F(\theta)$ is
algebraic over $F$, it follows that $\lambda$ is algebraic over $F$
and hence root of a non-zero polynomial $g(t)$ in $F[X]$. By applying
the automorphism $\sigma_i$ we conclude that $\sigma_i(\lambda)$ is
also root of $g(X)$ and hence $\sigma_i(\lambda)$, $1\leqq i\leqq n$,
is algebraic over $F$.

Now consider (***) as a system of linear equations with $\ell_0,
\dots, \ell_{n-1}$ as unknowns. The determinant of coefficient matrix
is the Vandermonde determinant $\prod_{j>i}(\theta_j-\theta_i)$, and
hence nonzero since the $\theta_i$'s are distinct elements of $M'$.
Consequently each $\ell_i$ belongs to $M'$.  Therefore $\ell_0,\dots,
\ell_{n-1}$ are algebraic over $F$ and thus lie in $F$.\hfill$\square$
 \enddemo

\remark{Remark 1.6} The separability assumption in Proposition 1.5
cannot be omitted, as the following example shows.  
Set $F_0=\Bbb F_p(a)$ where
$a$ is transcendental.  Then $a\in F_0$ and $\proot a \notin F_0$. Take
$s,t$ algebraically independent over $F_0$, and let 
 $$F=F_0(t), \qquad
M:=F(\proot t), \quad L:=F(s,\proot{g})\text{ where }
g=as^p+t.
 $$
 Then $M/F$ is inseparable. Clearly $\proot a\notin M$, but $(\proot
g-\proot t)/s=\proot a$ belongs to the compositum $LM$. Hence $M$ is
not algebraically closed relative to $LM$.  It is a cumbersome
computation with $p$'th powers and polynomials in two variables to
check that $F$ is algebraically closed relative to $L$. 

\ignore
 {\eightpoint\leftskip0,75truecm 
 We want to check that $F$ is algebraically closed relative to $L$. We
have $L=F(s)(\xi)$ with $\xi =\proot g$. Consider and element
$\alpha \in L$, say $\alpha =\varphi_0 +\varphi_1\xi
+\dots+\varphi_{p-1}\xi^{p-1}$ with $\varphi_i\in F(s)$.  Assume that
$ \alpha$ is algebraic over $F$.  We have to prove $\alpha \in
F$, that is, $\varphi _i=0$ for $i>0$ and $\varphi _0\in F$.

As $\alpha$ is algebraic over $F$, then so is 
 $$\alpha ^p=\varphi_0^p +\varphi_1^p
g+\cdots+\varphi_{p-1}^{\,p}g^{p-1}. \tag1
 $$
 But $\alpha ^p$  lies in   $F(s)$  and $F$ is
algebraically closed relative to $F(s)$; hence $\alpha ^p\in F$. 
As $F=F_0(t)$ we may write $\alpha ^p$ as a fraction of polynomials in
$F_0[t]$, and we may assume that the denominator is a $p$'th power.
Clearly, to prove the assertion we may multiply $\alpha$ with a
nonzero element of $F$.  Hence we may assume in (1) that $\alpha ^p$
is a polynomial in $F_0[t]$.

The coefficients $\varphi _i$ belong to $F(s)=F_0(t,s)$.  Write 
the $\varphi _i$ as fractions in their lowest terms, let $f\in
F_0[t,s]$ be the least common multiple of the denominators, and
multiply (1) with $f^p$.  The result is an equation
 $$f^p\alpha ^p=f_0^{\,p}+f_1^{\,p}g+\cdots+f_{p-1}^{\,p}g^{p-1},\tag2
 $$
 where each $f_i$ is a polynomial in $F_0[t,s]$ and no irreducible
factor of $f$ divides all the $f_i$. 

Clearly for any polynomial $h$ in $t$ we obtain, by grouping the
terms after
the powers of $t$, a unique expansion in the form
 $$h(t) =h ^{(0)}(t)+h ^{(1)}(t)t+\cdots h
^{(p-1)}(t)t^{p-1}, \tag3
 $$
 where each polynomial $h ^{(i)}(t)$ contains only powers of
$t^p$.  To obtain the expansion of $\alpha ^p$, substitute $s=0$ in
(2).  Indicate the substitution by a bar. Then $\bar
g=t$, and (2) yields  the expansion in $F_0[t]$,
 $$\bar f^p\alpha ^p=\bar f_0^{\,p}+\bar f_1^{\,p}t+\cdots+\bar
f_{p-1}^{\,p}t^{p-1},\tag4
 $$
 If $\bar f=0$, then by uniqueness of the expansion $\bar f_i=0$ for
$i=0,\dots,p-1$, that is, $s$ divides all the polynomials
$f,f_0,\dots,f_{p-1}$, contradicting the choice above.  Hence $\bar
f\ne 0$. Therefore, for the expansion of $h=\alpha ^p$ it
follows from (4) that $\bar f^ph^{(i)}=\bar f_i^p$.  Hence each
$h^{(i)}$ is a $p$'th power, say $h^{(i)}=k_i^{\,p}$ with
$k_i=F_0[t]$, and the expansion has the form,
 $$\alpha ^p= k_0^{\,p}+k_1^{\,p}t+\cdots k_{p-1}^{\,p}t^{p-1}.\tag 5
 $$
 Consider on the other hand the two sides of (2) as a polynomial
$\Phi$ in $F_0[s][t]$, and equate for the expansion the terms $\Phi
^{(0)}$.  The result is the equation,
  $$f^pk_0^p=f_0^{\,p}+f_1{\,p}(s^pa)+\cdots+f_{p-1}^{\,p}(s^pa)^{p-1},
 $$
 or
  $$0=(f_0-fk_0)^p+(f_1s)^pa+\cdots (f_{p-1}s^{p-1})^pa^{p-1}.\tag6
  $$
  Clearly $\proot a\notin F_0(s,t)$.  Hence the powers
$1,a,\dots,a^{p-1}$ are linearly independent over the subfield of
$p$'th powers of elements in $F_0(s,t)$.  Consequently the equation
(6) implies that $f_0=fk_0$, $f_1=\cdots=f_{p-1}=0$.  Hence $\varphi
_0=f_0/f=k_0\in F_0(t)$, and $\varphi _i=f_i/f=0$ for $i>0$, as
asserted, and the proof is complete.\hfill $\square$
 
}
\endignore
\endremark

\ignore
\remark{Drop this? Potential  examples}
 Question: 
Assume $L_1,L_2$ algebraically closed. Is the pair $L_1,L_2$ free over
$F=L_1\cap L_2$? Or linearly disjoint?

Potential example, far out: The field $F(s,t,u,v)$ ($F=\Bbb C$) where
$s^2+tu+v^2=0$  contains $L_1=F(s,t)$ and $L_2=F(u,v)$. Is 
$L_1\cap L_2=F$? Is $\widetilde L_1\cap \widetilde L_2=F$?
\endremark
\endignore

\subhead Section 2. Socles of field extensions
\endsubhead

We recall some basic concepts and facts from Galois theory.
  In the following $p$ always denotes a prime number. The cyclic group
of order $p$ is denoted $\Cycl_p$.  If $G$ is a group, by a {\it
$G$-extension of a field $F$} we mean a Galois extension of $F$ with
$G$ as Galois group.  

By an {\it elementary abelian $p$-extension\/}
of $F$ we mean a Galois extension whose Galois group is an elementary
abelian pro-$p$-group, that is, a product (finite or infinite) of
copies of $\Cycl_p$.

\definition {Definition 2.1} For an arbitrary field extension $L/F$
the {\it $p$-socle} $\Soc^p(L/F)$ is defined as the compositum of all
$\Cycl_p$-extensions of $F$ contained in $L$, or, equivalently, as the
largest elementary $p$-extension 
 of $F$ contained in $L$. (Hence the Galois group of  $\Soc^p(L/F)/F$
is a product of copies of $\Cycl_p$.)  
 \enddefinition

\definition{Definition 2.2} Let $G$ be a finite or a profinite group.
The {\it $p$-Frattini subgroup\/} $\Phi^p(G)$ is defined as the
intersection of all closed normal subgroups of index $p$ in $G$.
 The quotient $G/\Phi^p(G)$ is called the {\it $p$-Frattini
quotient\/} of $G$.

 More generally, for a subgroup $H$ of $G$ we denote by $\Phi ^p(G,H)$
the intersection
of all closed normal subgroups containing $H$ and of index $p$ in $G$.
Equivalently, $H\Phi ^p(G)$ is the intersection of the kernels of all
continuous homomorphisms $\varphi\:G\to \Cycl_p$ which are trivial on
$H$. Clearly, $\Phi(G,H)=H\Phi ^p(G)$  
 \enddefinition

Standard Galois theory immediately implies 
the following result.

\proclaim{Proposition 2.3} Let $M/F$ be a (not necessarily finite)
Galois extension with Galois group $G$. The $p$-socle of $M/F$ is the
fixed field of the $p$-Frattini subgroup $\Phi ^p(G)$, and the Galois
group of $\Soc^p(M/F)/F$ is the $p$-Frattini quotient $G/\Phi ^p(G)$.

More generally, if $L\subseteq M$ is the subfield fixed under a subgroup
$H\subseteq G$ then the $p$-socle of $L/F$ is the fixed field of the
subgroup $\Phi ^p(G,H)$.
\endproclaim

  \proclaim {Proposition 2.4} 
 Let $L/F$ and $M/F$ be two (not necessarily finite) Galois
extensions for which $L\cap M= F$.  Then the $p$-socle of the
compositum $LM/F$ is the compositum of the $p$-socles of $L/F$
and $M/F$.

 Assume more generally that $L/F$ and $M/F$ are separable extensions
for which $L'\cap M=F$, where $L'$ is the Galois
closure of $L$. Then the $p$-socle of the
compositum $LM/F$ is the compositum of the $p$-socles of $L/F$
and $M/F$.
 \endproclaim 

\demo{Proof of 2.4 in the Galois case} If $L/F$ and $M/F$ are finite
Galois extensions the assertion follows from classical Galois theory,
since obviously the $p$-Frattini quotient of a product is the product
of the $p$-Frattini quotients of the factors. For infinite Galois
extensions the assertion follows by writing the extensions as directed
unions of finite Galois extensions.  \hfill$\square$
 \enddemo

\demo{Proof of 2.4 in general} 
 The general result is not used in the proof of our main results, it
is more elaborate, and only included for completeness.  Also
for the general result we may assume that $L/F$ and
$M/F$ are finite separable extensions. Let $G$ be the Galois group of
a Galois extension of $F$ containing $LM$. Then the extension contains
$L'M$, and the subfields $L\subseteq L'$, and $M$, correspond to
subgroups, say, $N\supseteq N'$ and $H$ of $G$.  Clearly $LM$
corresponds to subgroup $N\cap H$. The subset $N'H\subseteq G$
is a subgroup since $N'$ is normal; it corresponds to the
intersection $L'\cap M$ which equals $F$ by assumption. Hence
$N'H=G$. In particular, $NH=G$, that is, any $g\in G$ has a
factorization,
 $$g=nh\text{ where }n\in N, \, h\in H.
 $$
 We use the factorizations to define, for any homomorphism $\varphi
\:G\to \Cycl_p$ which is trivial on $N\cap H$, two maps
$\varphi_N,\varphi_H \:G\to \Cycl_p$ such that
 $$\varphi _N(g)=\varphi (n),\ \varphi_H(g)=\varphi(h)\quad 
\text{when
}g=nh,\ n\in N,h\in H.
 $$
 First, the maps $\varphi _N,\varphi_H$ are well defined since
$\varphi $ is a homomorphism trivial on $N\cap H$.  Clearly $\varphi
(g)=\varphi _H(g)\varphi _N(g)$. Next, it is obvious that $\varphi_N$
equals $1$ on $H$ and $\varphi _H$ equals $1$ on $N$.  We show at
the end of the proof that the two maps $\varphi_N,\varphi _H\:G\to \Cycl_p$
are group homomorphisms.  

Now the equation of $p$-socles is equivalent to the following equation
of relative $p$-Frattini subgroups:
 $$\Phi ^p(G,\,N{\cap}H)=
\Phi^p (G,N)\cap \Phi ^p(G,H).\tag*
 $$
 As $\Phi^p(G,N)$ is the intersection of the kernels of all
homomorphisms $\psi \:G\to \Cycl_p$ trivial on $N$, it is obvious
that the inclusion ``$\subseteq$'' holds in (*).

To prove the reverse inclusion ``$\supseteq$'', let $g$ be an element
on the right side of (*).  We have to prove that $g$ belongs to the
left side, that is, $\varphi (g)=1$ for every homomorphism
$\varphi\:G\to \Cycl_p$ which is trivial on $N\cap H$. Use the above
decomposition $\varphi =\varphi _N\varphi _H$ where $\varphi _N$ is
trivial on $H$ and $\varphi _H$ is trivial on $N$. As shown below, the
maps are homomorphisms $\varphi_N,\varphi _H\:G\to\Cycl_p$. Therefore,
since $g$ belongs to the right hand side we have $\varphi _N(g)=1$ and
$\varphi _H(g)=1$.  Consequently $\varphi (g)=\varphi _N(g)\varphi
_H(g)=1$, as asserted.

It remains to prove that the maps $\varphi _N$ and $\varphi _H$ are
group homomorphisms.  So let $g_1,g_2\in G$, and write $g_1=n_1h_1$,
$g_2=n_2h_2$ with $n_1,n_2\in N$ and $h_1,h_2\in H$. As noticed in the
very beginning of the proof we have even $G=N'H$; in particular we may
assume that $n_2\in N'$. Now use the equation,
 $$g_1g_2=n_1h_1n_2h_2=(n_1h_1n_2h_1^{-1})(h_1h_2).
 $$
 Here $h_1h_2\in H$, and $n_1h_1n_2h_1^{-1}\in N$, since $\widetilde
N$ is normal and $N'\subseteq N$. Hence
 $$\varphi _H(g_1g_2)=\varphi (h_1h_2)=\varphi(h_1)\varphi(h_2)=\varphi
_H(g_1)\varphi _H(g_2).
 $$
 Similarly, since $\varphi _N$ takes values in a commutative group, it
follows that $\varphi _N$ is a homomorphism.
 \hfill $\square$
\enddemo

\example{Example 2.5} 
  Clearly the assumption $L'\cap M=F$ in 2.4 is stronger than
assuming that $L/F$ and $M/F$ are linearly disjoint. The stronger
assumption is only used at the very end of the proof to justify that
the maps $\varphi _N$ and $\varphi _H$ are homomorphisms.  The authors
haven't been able to decide whether the conclusion in 2.4 holds under
the simpler assumption that $L/F$ and $M/F$ are linearly disjoint.  In
the following example the two extensions are not linearly disjoint.

Let $F:=\Bbb Q(\varepsilon _p)$, where $\varepsilon _p$ is a primitive
$p$'th root of unity.  Choose a prime $q$ such that $q\equiv 1\pmod
p$, and let $L_1:=F(\root q\of 2)$ and $L_2=F(\varepsilon_q\root q\of
2)$.  Clearly $L_1\cap L_2=F$.  The extensions $L_1/F$ and $L_2/F$
have degree $q$, and hence their $p$-socles are trivial:
$\Soc^p(L_1/F)=\Soc^p(L_2/F)=F$. However, the compositum $L_1L_2$
contains the field $F(\varepsilon _q)$ which is a
$\Cycl_{q-1}$-extension of $F$. So there is a unique
$\Cycl_p$-extension of $F$ contained in $F(\varepsilon _q)$.  In
particular, the $p$-socle of $L_1L_2/F$ is nontrivial. 

However, the $p$-socle of $L_1L_2/F$ is nontrivial.  Indeed, the
compositum $L_1L_2$ contains the field $F(\varepsilon _q)$ which is a
$\Cycl_{q-1}$-extension of $F$. So there is a unique
$\Cycl_p$-extension of $F$ contained in $F(\varepsilon _q)$.  
 \endexample

\proclaim{Proposition 2.6} Assume that $E/F$ is contained in 
a Galois $p$-extension $N/F$.  If $\Soc^p(E/F)=F$, then $E=F$.
 \endproclaim
\demo{Proof} The assertion translates into the well-known 
property of $p$-Frattini subgroups of a $p$-group: If $H\subseteq G$
is an inclusion of $p$-groups, such that $H\Phi ^p(G)=G$ then $H=G$. 
 \hfill$\square$
 \enddemo

 \proclaim{Corollary 2.7} Let $N/F$ be a Galois $p$-extension and
$K/F$ an arbitrary extension. If $\Soc^p(N/F)\cap K=F$ then
$N\cap K=F$.
 \endproclaim
 \demo{Proof} Use the previous Proposition with $E:=N\cap K$. 
 \hfill$\square$
 \enddemo

\remark{Remarks on Kummer extensions and Artin--Schreier extensions
2.8}
 We fix a field $F$ and distinguish between the case, where the
characteristic of $F$ is $\neq p$ and the case where the
characteristic of $F$ is $p$.  

In the first case we will assume that $F$ contains the $p$-th roots of
unity.  In this situation any $\Cycl_p$-extension is a Kummer
extension $F(\root p \of a)$, $a\in F^*$.  Here $\root p \of a$, being
a root of $X^p -a$, is only determined up to factor $\zeta$ where
$\zeta$ is a $p$-th root of unity, but since $F$ contains $\zeta$ the
corresponding root fields coincide.

 For a subset $A\subseteq F^*$ we denote by $F(\root p\of A)$ the
extension generated by all elements $\root p\of a$ for $a\in A$.  In
particular, for an extension $L/F$ we have $\Soc^p(L/F)=F(\root p\of
A)$ where $A=(L^*)^p\cap F$. We shall need the following observations:
If $a_1,\dots,a_n$ are elements of $F^*$, then the extension $F(\root
p\of {a_1},\dots, \root p\of {a_n})/F$ has Galois group $\Cycl_p^n$ if
and only if the classes of the $a_i$ in the (multiplicative) group
$F^*/(F^*)^p$ are linearly independent over $\Bbb F_p$, 
 that is, if the relation,
 $$a_1^{\,\nu_1}\cdots a_n^{\,\nu_n}\in (F^*)^p, 
 $$
 implies that $\nu_i\equiv 0\pmod p$, $1\le i\le n$. 
  Also, if $A$
is a subset of $F^*$, and $b\in F^*$, then $\root p\of b\in F(\root
p\of A)$ if and only if the class of $b$ in group $F^*/(F^*)^p$
belongs to the subgroup generated by the classes of the elements of $A$.

\medskip
 With a similar notion in characteristic $p$ we have the following
observations: Any  $\Cycl_p$-extension of $F$
is an Artin--Schreier extension $F(\wp^{-1}a)$.  Here
$\wp(X)=X^p-X$ and
$\wp^{-1}a$, denoting a root of $X^p-X-a$, is only 
determined up to addition of a constant $c\in \Bbb F_p$. but the root
fields coincide.

For an extension $L/F$ we have $\Soc^p(L/F)=F(\wp^{-1}A )$ where
$A=\wp(L)\cap F$.  If $a_1,\dots,a_n$ are elements of $F$, then the
extension $F(\wpinv {a_1},\dots, \wpinv {a_n})/F$ has Galois group
$\Cycl_p^n$ if and only if the classes of the $a_i$ in the (additive)
group $F/\wp(F)$ are linearly independent over $\Bbb F_p$.  Also, if
$A$ is a subset of $F$, and $b\in F$, then $\wpinv b\in F(\wpinv A)$
if and only if the class of $b$ in group $F/\wp(F)$ belongs to the
subgroup generated by the classes of the elements of $A$. 
 \endremark

\remark{Remarks on maximal abelian extensions 2.9}
 Denote by $G_F$ the Galois group of the separable closure of $F$, and
by $G^{\operatorname{ab}}_F$ the Galois group of the maximal abelian
extension of $F$.  In addition, denote by $G_F(p)$ and
$G^{\operatorname{ab}}_F(p)$ the maximal pro-$p$ quotient groups.
Thus $G^{\operatorname{ab}}_F(p)$ is the Galois group of the maximal
abelian $p$-extension of $F$. 

 Assume the $F$ contains all $n$'th roots of unity for exponents $n$
not divisible by the characteristic of $F$. Then the group
$G^{\operatorname{ab}}_F(p)$ is a free pro-abelian $p$-group,
determined by the formula:
 $$G^{\operatorname{ab}}_F(p)=\hat Z_p^{I_p}\text{ where
}|I_p|=\cases \operatorname{rank}_{\FF_p} F^*/(F^*)^p&
\text{if }p\ne \operatorname{Char}(F),\\
 \operatorname{rank}_{\FF_p} F/\wp (F)&
\text{if }p= \operatorname{Char}(F),
\endcases\tag*
 $$
 where $\hat \ZZ_p$ is the (additive) group of $p$-adic integers.
Indeed, assume first that $p\ne \operatorname{Char}(F)$.  Then Kummer
theory may be combined with a theorem of  Capelli \cite{\Capelli}, see
also \cite{\Lorenz}. 
Accordingly, if $a\in K^*$ and $a\notin (K^*)^p$, then the $i$'th
field in the tower,
 $$F\subset F(\root p\of a)\subset \cdots \subset F(\root p^i\of
a) \subset\cdots\,,
 $$
 is a cyclic extension of $F$, and it is of degree $p^i$ by Capelli.
Consequently, the union of the fields in the tower is a $\hat
\ZZ_p$-extension of $F$.  It follows easily that the Galois group
$G^{\operatorname{ab}}_F(p)$ is the free pro-abelian $p$-group of rank
equal to the rank of $F^*/(F^*)^p$ as a vector space over $\FF_p$,
that is, the formula in the first case holds.

 The formula in the second case follows from the theorem of Witt:
The Galois group $G_F(p)$ of the maximal $p$-extension of $F$ is the
free pro-$p$-group of rank equal to the rank of $F/\wp (F)$.
Hence its maximal abelian quotient is a free pro-abelian group of the
same rank. 

\smallskip
 It follows from the formula that the Galois group
$G_F^{\operatorname{ab}}$ of the maximal abelian extension of $F$ is
the pro-abelian group,
  $$G^{\operatorname{ab}}_F=\prod_p \ZZ_p^{I_p},
 $$
where the product is over all primes $p$ and the cardinality of $I_p$
is given by (*).

\endremark

\subhead
Section 3.  Constructions of elementary abelian  $p$-extensions
\endsubhead

In this section we fix an algebraically closed field $C$, and we let
$T$ and $U$ be two nonempty subsets of a set
of independent variables over $C$.  We work inside a large
algebraically closed field containing the rational function field
$C(T,U)$, and consider the separable closures $\Coverline(T)$ and
$\Coverline(U)$, where separable closure has been indicated with a 
tilde. We show that the compositum $\Coverline(T)\Coverline(U)$ of
$\Coverline(T)$ and $\Coverline(U)$ is not separably closed and
that any finitely generated pro-nilpotent group and that any
elementary abelian $p$-group of rank at most the maximum of $|C|$,
$|T|$, and $|U|$, is realizable as Galois group over the
compositum.  We proceed in a series of lemmas.

\ignore 
 \csname proclaim\endcsname{Lemma 3.1} 
 The intersection $\Coverline(T)\,\cap\, \Coverline(U)$ equals $C$. 
 \endproclaim

Indeed, let $\alpha\in \Coverline(T)\,\cap\, \Coverline(U)$.  By
Proposition 1.4, $C(T)$ is algebraically closed relative to
$C(T)\Coverline(U)$. Hence $\alpha \in C(T)$.  Similarly $\alpha \in
C(U)$.  Consequently $\alpha \in C$. 

\csname proclaim\endcsname{Lemma 3.2}
The intersection $\Lambda
:=C(T)\Coverline(U)\,\cap\, \Coverline(T)C(U)$ equals $C(T,U)$. 
\endproclaim

Indeed, considering the diagram,
 $$\matrix
\Coverline(T) & \subseteq &\Coverline(T)C(U)\\
\cupeq&&\cupeq\\
C(T) & \subseteq & C(T,U).
\endmatrix\tag\dag
 $$
 the Translation Theorem, see \cite{\Lorenz} implies 
that $\Lambda =(\Lambda \cap \Coverline(T)\,) C(U)$. 

Clearly $\Lambda\supseteq C(T,U)$. If the inclusion were strict there
would be an element $\beta\in \Lambda\cap \Coverline(T)$ such that
$\beta \notin C(T)$. 
On the other hand, Proposition 1.4 implies that $C(T)$ is
algebraically closed relative to $C(T)\Coverline(U)$.  Since $\beta \in
C(T)\Coverline(U)$ we get a contradiction since $\beta \in
\Coverline(T)\setminus C(T)$. Thus $\Lambda =C(T,U)$. 

\csname proclaim\endcsname{Lemma 3.3}  For a given prime $p$ the
$p$-socle of the compositum is given by the formula,
 $$\Soc^p\bigl(\Coverline(T)\,\,\Coverline(U)/C(T,U)\bigr)
=\bigl(\Soc^p\bigl(\Coverline(T)/C(T)\bigr)\bigr)\bigl(
\bigl(\Soc^p\bigl(\Coverline(U)/bigr)\bigr).
 \tag\ddag
 $$
 \endproclaim

Indeed, consider first the $p$-socles of $\Coverline(T)C(U)/C(T,U)$
and $C(T)\Coverline(U)/C(T,U)$. By the Translation Theorem applied to
the diagram (\dag) we get the first of the following two equations:
 $$\gather 
\Soc^p_{C(T,U)}(\Coverline(T)C(U))
 =(\Soc^p_{C(T)}(\Coverline(T))\,C(U),\\
\Soc^p_{C(T,U)}(C(T)\Coverline(U))=C(T)(\Soc^p_{C(U)}(\Coverline(U)).
\endgather
 $$
 The second equation is obtained similarly. From the two equations and
Proposition 2.4 we conclude (\ddag).
 \endignore

 \proclaim{Lemma 3.1} 
 The intersection $\Coverline(T)\,\cap\, C(T)\Coverline(U)$
equals $C(T)$. 
 \endproclaim
\demo{Proof}
The field $C$ is algebraically closed, and in particular algebraically
closed relative $\Coverline(U)$. By Proposition 1.4, $C(T)$ is
algebraically closed relative to $C(T)\Coverline(U)$. In particular,
if an element of $C(T)\Coverline(U)$ is separably algebraic over
$C(T)$ then it belongs to $C(T)$. 
\hfill$\square$
 \enddemo

 \proclaim{Lemma 3.2} 
 The intersection $\Coverline(T)C(U)\,\cap\, C(T)\Coverline(U)$
equals $C(T,U)$. 
 \endproclaim
\demo{Proof}
Consider the following diagram of fields and inclusions:
 $$\matrix
\Coverline(T) & \subseteq &\Coverline(T)C(U)&\subseteq&
                   \Coverline(T)\Coverline(U)\\
\cupeq&&\cupeq&&\cupeq\\
C(T) & \subseteq & C(T,U) & \subseteq& C(T)\Coverline(U).
\endmatrix\tag\dag
 $$
 To facilitate the notation, we let $L_0:=C(T)$, 
$M_0:=C(U)$, $F:=C(T,U)$, $M:=C(T)\Coverline(U)$ and we let
$L:=\Coverline(T)C(U)=\tilde L_0M_0$.  Then the fields of the above
diagram are the following:
 $$\matrix
\tilde L_0 &\subseteq &L &\subseteq &LM\\
\cupeq&&\cupeq&&\cupeq\\
L_0 & \subseteq & F & \subseteq & M
\endmatrix,\tag\ddag
 $$
 and $LM=\tilde L_0\tilde M_0$. Consider the
vertical extensions.  The first is Galois; its Galois group
$G_0:=\Gal(\tilde L_0/L_0)$ is the absolute Galois group of $L_0$.
Therefore, by the Translation Theorem, see \cite{\Lorenz}, the next two
vertical extensions are also Galois, and for the Galois groups
$G:=\Gal (L/F)$ and $G':=\Gal(LM/M)$ we have injections,
 $$G_0 \hookleftarrow G \hookleftarrow G'.
 $$
 The inclusion $G'\to G_0$ is surjective if and only if $\tilde
L_0\cap M=L_0$. Hence, by Lemma 3.1, the inclusion $G'\to G_0$ is an
isomorphism.  Therefore, so is the inclusion $G'\to G$. Hence, again
by the Translation Theorem, we have the equation $L\cap M=F$, which is
the asserted equality.
\hfill$\square$
 \enddemo 

 \proclaim{Lemma 3.3}  For the given
prime $p$ we have the following
equation of $p$-socles:
 $$\Soc^p\bigl(\Coverline(T)\Coverline(U))\,/\,C(T,U)\bigr)
=\Soc^p\bigl(\Coverline(T)\,/\,C(T)\bigr)\,\,
  \Soc^p\bigl(\Coverline(U)\,/\,C(U)\bigr).
 $$
 \endproclaim
 \demo{Proof} We use the same notation as in the previous proof, and
refer to the diagram (\ddag). By the previous Lemma, the Galois group
$G=\Gal(L/F)$ equals $G_0=\Gal (\tilde L_0/L_0)$, the absolute Galois
group of $L_0$. So the two groups have the same $p$-Frattini
quotients. Consequently, for the $p$-socle we obtain the first of the
following equations,
 $$\aligned
\Soc^p (L/F)&=\Soc^p(\tilde L_0/L_0)F,\\
\Soc^p(M/F)&=\Soc^p(\tilde M_0/M_0)F .\endaligned\tag*
 $$
 The second equation follows from the
symmetry in $T$ and $U$.

Finally, by Proposition 2.4, the $p$-socle of $LM/F$ is the compositum
of the $p$-socles of $L/F$ and $M/F$, or, by the equations (*),
 $$\Soc^p(LM/F)=\Soc^p(\tilde L_0/L_0) \Soc^p(\tilde M_0/M_0).
 $$
 The latter equation is the asserted equality of socles.
 \hfill$\square$\enddemo

 \proclaim{Lemma 3.4} Assume that $p\ne \operatorname{Char}C$.
Fix $t\in T$ and $u\in U$, and consider the field $C(t+u)$. Take any
$n$ distinct elements $c_1,\dots,c_n\in C$. Then the following
extension:
 $$C(t+u)\bigl(\proot{t+u+c_1},\dots,\proot{t+u+c_n}\,\bigr)
 $$
 is an elementary abelian  $p$-extension of $C(t+u)$ with Galois group
$\Cycl_p^{\,n}$, and linearly disjoint with
$\Coverline(T)\Coverline(U)$ over $C(t+u)$. 
 \endproclaim
\demo{Proof}
 For simplicity, set $q_i:=t+u+c_i$, $i=1,\dots,n$. It suffices to
prove that the $n$ roots $\proot{q_1},\dots,\proot{q_n}$ generate over
$\Coverline(T)\Coverline(U)$ an elementary abelian $p$-extension with
group $\Cycl_p^{\,n}$. The polynomials $q_i$ belong to the field
$C(T,U)$, and by the remarks on Kummer theory in 2.8, it suffices to
prove that the classes of the $q_i$ modulo the multiplicative subgroup
of $C(T,U)^*$ describing the $p$-socle of $\Coverline(T)\Coverline(U)$
over $C(T,U)$ are linearly independent over $\Bbb F_p$.  The $p$-socle
is described in Equation (\ddag) in Lemma 3.3, and the socles on the
right hand side of the equation are the maximal elementary abelian
$p$-extensions of $C(T)$ and $C(U)$.  Hence it suffices to show for
any equation of the following form in $C(T,U)^*$:
 $$q_1^{\nu_1}\cdots q_n^{\nu_n}
=\varphi\psi\alpha^p\qquad(\text{with } \nu_i\in \Bbb
Z),\tag**
 $$
 where $\varphi \in C(T)^*$, $\psi \in C(U)^*$ and $\alpha \in
C(T,U)^*$, that $\nu_i\equiv 0\pmod p$ for $i=1,\dots,n$.  To prove
it, fix $i$, and let $v\:C(T,U)^*\to \Bbb Z$ be the $q_i$-adic
valuation. Then $v(q_j)=0$ for $j\ne i$, and $v(q_i)=1$.  So the value
of $v$ on the left of $(**)$ is $\nu_i$.  On the right, the value of
$v$ at $\varphi$ and $\psi$ is zero, since $q_i$ does not belong to
$C(T)$ or to $C(U)$.  Hence the value on the right side is a multiple
of $p$.  Therefore $\nu_i\equiv 0\pmod p$.
 \hfill$\square$
 \enddemo

\ignore
 \csname proclaim\endcsname{Lemma 3.4} Any elementary abelian
pro-$p$-group {\rm(}for $p\ne \operatorname{Char}C${\rm)} of rank at
most the maximum of $|C|$, $|T|$, $|U|$, can be realized as a Galois
group over the compositum $\Coverline(T)\Coverline(U)$. 
 \endproclaim

 Indeed, since $C$ is algebraically closed, all $p$'th roots of unity
lie in $C$. Take any $t\in T$, $u\in U$, and $n$ distinct elements
$a_1,\dots,a_n\in C$, and consider the $n$ elements of $C(T,U)$:
 $$q_i :=t+u+a_i, \quad i=1,\dots, n
 $$
 It is easy to see that the $p$'th roots $\root p\of {q_1 }, \dots,\root
p\of {q_n }$ generate a $p$-extension of
$\Coverline(T)\cdot\Coverline(U)$ with Galois group $\Cycl_p^{\,n}$,

{\eightpoint\leftskip0,75truecm
 In fact, by Proposition 2.4 and the remarks in 2.8 it suffices to
prove for any equation of the following form in $C(T,U)^*$ 
 $$q_1^{\nu_1}\cdots q_n^{\nu_n}
=\varphi\psi\alpha^p\qquad(\text{with } \nu_i\in \Bbb
Z),\tag***
 $$
 where $\varphi \in C(T)^*$, $\psi \in C(U)^*$ and $\alpha \in
C(T,U)^*$, that $\nu_i\equiv 0\pmod p$ for $i=1,\dots,n$.  To prove
it, fix $i$, and let $v\:C(T,U)^*\to \Bbb Z$ be the $q_i$-adic
valuation. Then $v(q_j)=0$ for $j\ne i$, and $v(q_i)=1$.  So the value
of $v$ on the left of $(***)$ is $\nu_i$.  On the right, the value of
$v$ at $\varphi$ and $\psi$ is zero, since $q_i$ does not belong to
$C(T)$ or to $C(U)$.  Hence the value on the right side is a multiple
of $p$.  Therefore $\nu_i\equiv 0\pmod p$.
  \par}
\endignore

\proclaim{Proposition 3.5}
 Any elementary abelian pro-$p$-group of rank at
most the maximum of $|C|$, $|T|$, $|U|$, can be realized
as a Galois group over the compositum
$\Coverline(T)\Coverline(U)$. 
 \endproclaim
 \demo{Proof}
 Assume first that $p\ne\operatorname{Char}C$. Then Lemma 3.4 applies.
Varying $n$ it follows first that we may realize $\Cycl_p^{|C|}$ over
the compositum, and varying similarly $t\in T$ and $u\in U$ we may
realize all elementary abelian pro-$p$-groups of the asserted rank. 

The case $p=\operatorname{Char}C$ remains. Only the first part of the
proof of Lemma 3.4 needs to be changed.  Take any $t\in T$, $u\in U$,
and set $q:=t+u$.  Let $c_1,\dots,c_n$ be any $n$ elements of $C$
linearly independent over the prime field $\Bbb F_p$.  We prove that
the $n$ roots $\wp^{-1}(c_1/q),\dots,\wp^{-1}(c_n/q)$ generate an
elementary abelian  $p$-extension of $\Coverline(T)\Coverline(U)$ with
Galois group $\Cycl_p^{\,n}$.  Using the description of the $p$-socle
in Lemma 3.4, and the remarks on Artin--Schreier theory in 2.8 is
suffices to show for any
equation of the following form in $C(T,U)$:
 $$\nu_1c_1/q+\cdots+ \nu_nc_n/q =\varphi+\psi
+\alpha ^p-\alpha\qquad (\text{with }\nu_i\in \Bbb
Z),
 $$
 where $\varphi \in C(T)$, $\psi \in C(U)$ and $\alpha \in C(T,U)$,
that $\nu_i\equiv 0\pmod p$ for $i=1,\dots,n$. It suffices to prove
that the equation implies that the numerator
$c=\nu_1c_1+\cdots+\nu_nc_n$ on the left side vanishes.  Assume that
$c\ne 0$.  Then the value of the $q$-adic valuation on the left side
equals $-1$.  On the right side the value is zero or $+\infty$ at
$\varphi$ and at $\psi$.  Hence a contradiction is obtained since no
discrete valuation can take the value $-1$ at an element of the form
$\alpha^p-\alpha=\alpha(\alpha^{p-1}-1)$.
 \hfill$\square$
 \enddemo

\ignore

\remark{Remark 3.6}
  Let $N$ be the compositum $N:=\Coverline(T)\Coverline(U)$. Fix $t\in
T$ and $u\in U$, and set $x:=t+u$.  Then the subfield $C(x)\subseteq
N$ is the field of rational functions in one variable over $C$. The
first part of argument in 3.4 implies the following: Any elementary
abelian $p$-extention $S/C(x)$ is linearly disjoint with $N/C(x)$.
Equivalently, if $S/C(x)$ is an elementary abelian $p$-extension with
Galois group $G$, then $SN/N$ is a $G$-extension.

Indeed, it suffices to consider a finite extension $S/C(x)$. 
As $C$ is algebraically closed, the polynomials $x+a$ for $a\in C$
generate the group $C(x)^*$ modulo $C^*$. Therefore $S$ is contained
in the extension $C(x)(\proot {q_1},\dots,\proot{q_n})$ for suitable
$n$. By the argument of 3.4, the latter extension is linearly disjoint
with $N$. Therefore, so is the former.

 What a pity:  Witt implies that rest of this remark may be dropped: 
 It is possible that an almost similar result holds when
$p=\operatorname{Char}(C)$. Indeed, by decomposition of rational
functions, an $\Bbb F_p$-basis of the additive group $C(x)$ is
obtained as follows: Let $(c_j)_{j\in J}$ be an $\Bbb F_p$-basis for
$C$.  Then the union of the two series of elements,
 $$\gamma _{j\nu}:=c_jx^\nu, \ j\in J,\nu \ge 0\ \text{ and }\
\gamma_{j,\nu,a}=\frac{c_j}{(x+a)^\nu},\ j\in J,\nu\ge 1, a\in C.
 $$
 form a basis. It is easy to see that modulo the image of $\wp$ a
basis is formed by taking the elements of the series for $1\le \nu \le
p-1$.  It follows from arguments similar to those used in 3.5 that the
elements of the second series remain independent in $N/\wp(N)$.

But this is not true for elements of the first series.  For instance,
$\wp^{-1}cx=\wp^{-1}ct+\wp^{-1}cu$ belongs to
$\Coverline(T)\Coverline(U)$. 

QUESTION.  How do we apply this?  Characterize elementary abelian
$p$-extensions of $C(t)$ contained in $C(t)(\wp^{-1}(\gamma_{j,\nu,a})$.  Can we
apply the D-H-P theorem in this case.
 \endremark

\endignore

\subhead
Section 4.  Realization of Galois groups over the compositum
\endsubhead

In this section we consider two algebraically closed fields $L$ and
$M$ contained in some common larger algebraically closed field.  We
assume that none of the two fields is contained in the other, and we
assume that $L$ and $M$ are linearly disjoint over their intersection
$C=L\cap M$. Clearly the intersection $C$ is algebraically closed.  We
show that the compositum $LM$ is not separably closed and, moreover,
we show that any pro-nilpotent group with a generating set of
cardinality at most $|C|$ is realizable as a Galois group over $LM$.
In addition, any elementary abelian $p$-group with a generating set of
cardinality at most the maximum of $|L |$ and $|M|$ is realizable over
$LM$. 

Since $C$ is perfect, we can chose separating transcendency bases $T$
for $L/C$ and $U$ for $M/C$, see \cite{\ZS, Theorem 31, p.~105}. Both
are non-empty since $C$ is strictly contained in $L$ and $M$.
Moreover, since $L/C$ and $M/C$ are linearly disjoint, it follows from
the observations in Remark 1.2 that the union $T\cup U$ is
algebraically independent over $C$.  So the setup of Section 3
applies: $L=\Coverline(T)$ and $M=\Coverline(U)$, except that the two
fields were denoted $\tilde L_0$ and $\tilde M_0$ in Lemma 3.2.
Clearly the maximum of $|C|$, $|T|$, and $|U|$ is equal to the maximum
of $|L|$ and $|M|$,  and hence equal to $|LM|$.

\proclaim{Theorem 4.1} Let $L$ and $M$ be algebraically closed fields,
of which none is contained in the other, and linearly disjoint over
their intersection $C:=L\cap M$. Let $N:=LM$  be their
compositum.  Then the Galois group $G^{\operatorname{ab}}_N$ of the
maximal abelian extension of $N$ is the free pro-abelian group of rank
$|N|$, that is, $G^{\operatorname{ab}}_N\simeq \widehat Z^N$.  
\endproclaim
\demo{Proof}
 All roots of unity belong to $N$ since $N$ contains an algebraically
closed field.  Hence it follows from the remarks in 2.9 that
 $G^{\operatorname{ab}}_N(p)$ is a free pro-abelian $p$ group and that
its rank is given by the formula (*) is 2.9.  In the notation above
$N=\Coverline(T)\Coverline(U)$.  Hence it follows from 3.5 that the
rank is at least the maximum of $|C|$, $|T|$, and $|U|$, and hence at
least equal to $|N|$.  Consequently, the rank is equal to $|N|$, that
is, $G^{\operatorname{ab}}_N(p)=\hat \ZZ_p^N$.  Therefore, by formula
(\ddag) in 2.9, $G^{\operatorname{ab}}_N=\widehat \ZZ^N$, and
the theorem has been proved.  
 \hfill$\square$
 \enddemo

 \proclaim{Theorem 4.2} 
 In the setup of Theorem\/ {\rm (4.1)} the free pro-$p$-group 
of rank $|C|$ is realizable 
as Galois group over the compositum
$N=LM$.  As a consequence, the free pro-nilpotent group of rank $|C|$ is
realizable. 
 \endproclaim
 \demo{Proof}
 Fix from the separating transcendency bases elements $t\in T$ and
$u\in U$, let $K$ be the maximal $p$-extension of $C(t+u)$, and denote
by $G(p)$ the Galois group of $K/C(t+u)$. 

Assume first that $p\ne \operatorname{Char} C$.  Let $S$ be the
$p$-socle of $K/C(t+u)$.  The polynomials $t+u+c$, $c\in C$, generate
the multiplicative group $C(t+u)^*$ (up to multiplication by a
constant in $C^*$), since $C$ is algebraically closed.  Hence, by
Kummer theory, $S$ is the maximal extension of the form in Lemma 3.4.
The latter extension is linearly disjoint with $N$ over $C(t+u)$ by
Lemma 3.4.  So the subfield $S$ is linearly disjoint with $N$ over
$C(t+u)$.  Finally, by Proposition 2.6, $K$ is linearly disjoint with
$N$ over $C(t+u)$.  Consequently $\Gal(KN/N)=\Gal(K/C(t+u))=G(p)$.
Finally, it follows from the Douady--Harbater--Pop theorem on the
absolute Galois group of the function field $C(t+u)$ (see
\cite{\Douady},\cite{\Harbater},\cite{\Pop}, or Haran--Jarden
\cite{\HaranJarden} for an easily accessible proof) that 
 $G(p)$ is the free pro-$p$ group of rank $|C|$.

If $p=\operatorname{Char}(C)$, we use the Theorem of Witt \cite{\Witt}
(see also J.-P.~Serre \cite{\Serre, Corollaire 1, p.~91}) on the
maximal $p$-extension of a field of characteristic $p$.  Accordingly,
the maximal Galois $p$-extension of $N$ is a free pro-$p$-group
$G_N(p)$ of rank equal to the rank of $N/\wp(N)$.  It follows from
Proposition 3.5 that $|N/\wp(N)|\ge |C|$.  Consequently, the free
pro-$p$ group of rank $|C|$ is a quotient of $G_N(p)$, and hence
realizable over $N$.

Thus the first assertion has been proved for all $p$.  Clearly the
last assertion is a consequence since a pro-nilpotent group is the
over all primes $p$ of pro-$p$ groups, and hence realizable by the
compositum of the fields realizing the factors.
  \hfill$\square$
 \enddemo

 \remark{Note 4.3} Given the algebraically closed fields $L$ and $M$
none of which is contained in the other so that
$C=L\cap M$ is properly contained in them both.  It is
part of the results in Section 3 that $L/C$ and $M/C$ are linearly
disjoint if and only if the union $T\cup U$ of the separating
transcendency bases is algebraically independent over $C$.  
The condition of being linearly disjoint over $C$ is clearly
equivalent to the condition on the separating transcendency bases $T$
and $U$ that their union $T\cup U$ is algebraically independent over
$C$. 

The authors have not been able to decide if the condition
of linear disjointness is always
satisfied.  It  is clearly satisfied if one of $T$ and $U$ is a
singleton.  As a consequence, the conclusions in 4.1 and
4.2 hold if $L$ or $M$
has transcendence degree $1$ over $C$. 
\endremark

 \enddocument